\documentclass{amsart}
\usepackage{amssymb, amsmath, amsfonts, amscd}

\input xypic

\theoremstyle{plain}
\newtheorem{theorem}{Theorem}[section]
\newtheorem{corollary}[theorem]{Corollary}
\newtheorem{lemma}[theorem]{Lemma}
\newtheorem{proposition}[theorem]{Proposition}

\theoremstyle{definition}
\newtheorem{definition}[theorem]{Definition}

\theoremstyle{remark}

\numberwithin{equation}{theorem}

\renewcommand{\O}{\mathcal{O} }

\renewcommand{\P}{\mathbf{P} }
\renewcommand{\Pr}{\mathcal{P} }

\newcommand{\mf}[1]{\mathfrak{#1} }

\newcommand{\Der}{\operatorname{Der} }
\newcommand{\Hom}{\operatorname{Hom} }

\newcommand{\End}{\operatorname{End} }
\newcommand{\Spec}{\operatorname{Spec} }

\renewcommand{\H}{\operatorname{H} }
\newcommand{\DR}{\operatorname{DR}}

\newcommand{\da}{\alpha}
\newcommand{\db}{\beta}
\newcommand{\dc}{\gamma}

\newcommand{\K}{\operatorname{K} }
\newcommand{\C}{\operatorname{C} }
\newcommand{\R}{\operatorname{R} }
\newcommand{\Gal}{\operatorname{Gal} }

\renewcommand{\K}{\operatorname{K} }
\renewcommand{\lg}{ \mathfrak{g}}

\newcommand{\dd}{\nabla \oplus \nabla'}
\newcommand{\ddd}{\nabla \oplus \nabla'_\delta}
\newcommand{\dde}{\nabla \oplus \nabla'_\eta}
\newcommand{\dtd}{\nabla \otimes \nabla'}
\newcommand{\dtdd}{\nabla \otimes \nabla' _\delta}
\newcommand{\dtde}{\nabla \otimes \nabla'_\eta}

\newcommand{\nd}{\nabla_\delta}
\renewcommand{\ne}{\nabla_\eta}

\begin{document}

\title{The Chern character for Lie-Rinehart algebras}
\author{Helge Maakestad }
\address{Institut de Mathematiques, Universite Paris VII}
\email{maakestad@math.jussieu.fr} 

\keywords{Lie-Rinehart algebra, connections, algebraic stacks, differential graded algebras, Grothendieck rings, Chern characters, 
de Rham cohomology, Lie-Rinehart cohomology, Jacobson Galois correspondence }
\subjclass{14C17, 19E15, 12.80, 14L15}
\date{Spring 2005}
\begin{abstract} We construct a Chern character for the situation 
\[ ch^\lg:\K_0(\lg)\rightarrow \H^*(\lg, A) \]
where $\lg$ is any Lie-Rinehart algebra, $A$ is any $k$-algebra of characteristic zero and $\H^*(\lg, A)$ is the
Lie-Rinehart cohomology of $\lg$. Our result generalize the classical Chern character 
\[ ch^A :\K_0(A)\rightarrow \H^*_{\DR}(A), \]
where $\K_0(A)$ is the $\K$-theory of $A$ and $\H^*_{\DR}(A)$ is the algebraic de Rham cohomology.
\end{abstract}

\maketitle

\tableofcontents

\section*{Introduction} 

Classical Galois theory setting up a one to one correspondence between intermediate field-extensions
of a Galois extension $E\subseteq F$ and subgroups of the
Galois group $\Gal (F/E)$ was generalized by 
N. Jacobson in \cite{jacobson} to give a Galois-correspondence for purely inseparable field-extensions $k\subseteq K$ of exponent one 
of a field $k$ of characteristic $p>0$. This is a one to one correspondence between intermediate fields 
and $p-K/k$-sub-Lie algebras of $\Der_k(K)$. 
Jacobsons $p-K/k$-Lie algebra is the characteristic $p$  version  of a
structure called a Lie-Rinehart algebra. 

For an arbitrary $k$-algebra $A$, there exists the notion of a
$(k,A)$-Lie-Rinehart algebra:
it is a $k$-Lie algebra
and $A$-module $\lg$ with a map of $k$-Lie algebras and $A$-modules
$\alpha :\lg \rightarrow \Der_k(A)$, i.e
a Lie-algebra acting on $A$ in terms of vector fields. There exists the notion of a $\lg$-connection $\nabla$ on an  $A$-module $W$:
this is an action 
\[ \nabla:\lg \rightarrow \End_k(W) \]
generalizing the notion of a covariant derivation. There exists a complex $\C^\bullet(\lg,W ,\nabla )$ - the Lie-Rinehart 
complex - generalizing
simultaneously the algebraic de Rham complex of $A$  and the Chevalley-Eilenberg complex of $\lg$.
The main result of this paper is the following (see Theorem \ref{theoremtwo}): 
There exists a ring homomorphism
\[ ch^\lg :\K_0(\lg)\rightarrow \H^*(\lg,A) \]
from the Grothendieck ring  $\K_0(\lg)$ 
to the cohomology ring $\H^*(\lg,A)$. Here $\K_0(\lg)$ is the Grothendieck ring of locally free $A$-modules with a $\lg$-connection and 
$\H^*(\lg, A)$ is the Lie-Rinehart cohomology of $A$ with respect to $\lg$. 
We prove furthermore in Theorem \ref{independence} that the Chern character from Theorem \ref{theoremtwo} 
is independent with respect to choice of $\lg$-connection. 
This generalizes the construction of the classical Chern character (see Corollary \ref{corollaryone}.) 
Note that J. Huebschmann has in \cite{huebsch0} considered a Chern-Weil construction in a similar situation, 
and it would be interesting to relate the construction in \cite{huebsch0}  to the construction in this note. 

The notion of a $(k,A)$-Lie-Rinehart algebra is closely related to the notion of a groupoid in schemes. 
One constructs from a groupoid in schemes a Lie-Rinehart algebra in the same way as one constructs
the Lie algebra from a group scheme. Much of the theory for group schemes and Lie algebras 
can be generalized to this new situation. 

Lie-Rinehart algebras appear in topology and knot theory: T. Kohno has in \cite{KOHNO} computed 
the Alexander polynomial of an irreducible plane curve $C$ in $\mathbf{C}^2$ using the 
logarithmic deRham complex $\Omega^\bullet_{\mathbf{C}^2}(*C)$ which is just the standard complex where
we let $\lg$ be the Lie algebra of derivations preserving the ideal of $C$ in $\mathbf{C}^2$. 

Groupoids and Lie-Rinehart algebras appear in the theory of motives: Let $\mathbf{T}$ be a Tannakian category over a field $F$ 
of characteristic zero, and let $\omega$ be a fiber functor over the algebraic closure $\overline{F}$ of $F$. 
Then $Aut^\otimes (\omega)$ is represented by a groupoid $S/S_0$ and there exists an equivalence of categories
\[ \mathbf{T}\cong \mathbf{Rep}(S/S_0) ,\]
(see \cite{milne}).

The paper is organized as follows: In the first section we define and sum up various general properties
of Lie-Rinehart algebras, connections and the Lie-Rinehart complex. In the second section we prove existence
of the Chern character. In the third section we prove that the Chern
character is independent with respect to choice of connection.

\section{Lie-Rinehart algebras, connections and the Lie-Rinehart complex}

In this section we introduce objects in the theory of modules on
Lie-Rinehart algebras and state some general facts on the following: Let $A$ be a commutative
ring over a field $k$. Let furthermore
$\mf{g}$ be an $(k,A)$-Lie-Rinehart algebra and let $(W,\nabla )$ be a $\mf{g}$-module.
We introduce the Lie-Rinehart complex $\C^\bullet (\lg ,W, \nabla)$.
If $\nabla$ is flat, $\C^\bullet(\mf{g},W, \nabla )$ is a DG-module, hence $\H^\bullet(\mf{g},W, \nabla)$ is a
graded left $\H^\bullet(\mf{g},A)$-module.

\begin{definition} Let $A$ be a commutative $k$-algebra where $k$ is a
  commutative ring. A
\emph{$(k,A)$-Lie-Rinehart algebra} on $A$ is a $k$-Lie algebra and an $A$-module
$\mathfrak{g}$ with a map $\alpha:\lg\rightarrow \Der_k(A)$ satisfying
the following properties:
\begin{align}
&\alpha(a\delta)=a\alpha(\delta)   \\
&\alpha([\delta, \eta])=[\alpha(\delta),\alpha(\eta)] \\
&[\delta,a\eta]=a[\delta,\eta]+\alpha(\delta)(a)\eta 
\end{align}
for all $a\in A$ and $\delta,\eta \in \lg$.
Let $W$ be an $A$-module. A $\lg$-\emph{connection} $\nabla$ on $W$,
is
an $A$-linear map $\nabla:\lg\rightarrow \End_k(W)$ which satisfies
the \emph{Leibniz-property}, i.e. 
\[ \nabla (\delta)(aw)=a\nabla (\delta)(w)+\alpha (\delta)(a)w \]
for all $a\in A$ and $w\in W$. The $\lg$-connection $\nabla$ is \emph{flat} if it is a map of Lie algebras.
If $\nabla$ is flat, we say that the pair $(W,\nabla)$ is a \emph{$\lg$-module}.  
\end{definition}
When it is clear from the context the notion Lie-Rinehart algebra will
be use instead of $(k,A)$-Lie-Rinehart algebra. 
A Lie-Rinehart algebra is also referred to as a a \emph{Lie-Cartan pair} or a 
\emph{foliation}.

\begin{definition} \label{CEcomplex} Let $A$ be a $k$-algebra, $\lg$ a
  Lie-Rinehart algebra
and $(W,\nabla )$ an $A$-module with a $\lg$-connection. Define a sequence of $A$-modules  
$\tilde{C}^\bullet (\lg, W,\nabla )$ and $k$-linear differentials $d^\bullet$  in the following way:
Let $\tilde{C}^p(\lg,W,\nabla )=\Hom_k(\wedge^p\lg,W)$ where $\wedge^p\lg$ is wedge product over $A$. Define differentials
\[ d^p:\tilde{C}^p(\lg,W,\nabla )\rightarrow \tilde{C}^{p+1}(\lg,W,\nabla ) \]
by
\begin{align}\label{differential} 
&(d^p\psi)(\delta_1\wedge \cdots \wedge \delta_{p+1})=
\sum_{i=1}^{p+1}(-1)^{i+1}\nabla_{\delta_i}\psi(\delta_1\wedge \cdots
\wedge \hat{\delta_i}\wedge \cdots \wedge \delta_{p+1})
\end{align} 
\[+ \sum_{1\leq i<j\leq p+1}(-1)^{i+j}\psi([\delta_i,\delta_j]\wedge
\cdots \wedge \hat{\delta_i} \wedge \cdots \wedge \hat{\delta_j}
\wedge \cdots \wedge \delta_{p+1}) .\]
Put $\tilde{C}^0=W$ and define $d^0(w)(\delta)=\nabla(\delta)(w)$. Let $R_\nabla=
d^1\circ d^0$ be the \emph{curvature} of the connection $\nabla$.
\end{definition}
Notice that $R_\nabla (\delta \wedge \eta)=[\nabla _\delta,\nabla _\eta]-\nabla_{[\delta,
\eta ]}$ hence $W$ is a $\lg$-module if and only if
the curvature is zero.  Note also:  if the connection $\nabla$ is flat and $A=k$, the sequence of modules and differentials 
defined in \ref{CEcomplex} is just the ordinary \emph{Chevalley-Eilenberg} complex of the representation $W$ 
for the $k$-Lie algebra $\lg$.

\begin{lemma}\label{complex} 
Let  $\lg$ be a Lie-Rinehart algebra and let $(W,\nabla)$ be a $\lg$-connection. 
Consider the sequence of modules from definition \ref{CEcomplex}, 
  $\tilde{C}^\bullet(\lg,W,\nabla)$. Then for all $p\geq 0$ the
  following holds: 
\[(d^{p+1}\circ d^p)(\delta_1\wedge \cdots \wedge \delta_{p+2}) \]
\[= \sum_{1\leq i < j \leq p+2}(-1)^{i+j+1}R_\nabla (\delta_i\wedge
  \delta_j)
(\delta_1\wedge \cdots \wedge \hat{\delta_i}\wedge \cdots \wedge \hat{\delta_j}
  \wedge \cdots \wedge \delta_{p+2}) .\]
Furthermore the maps $d^p$ induce maps
\[ d^p:\Hom_A(\wedge^p\lg,W)\rightarrow \Hom_A(\wedge^{p+1}\lg,W) \]
i.e $d^p\phi(aw)=ad^p\phi(x)$. 
\end{lemma}
\begin{proof} Standard fact. 
\end{proof}

\begin{definition}\label{standardcomplex} Define the
  \emph{Lie-Rinehart complex} $C^\bullet(\lg ,W,\nabla)$ as follows:
\[ C^p(\lg , W, \nabla)=\Hom_A(\wedge^p \lg, W) ,\]
with differentials
\[ d^p: \Hom_A(\wedge^p \lg, W)\rightarrow \Hom_A(\wedge^{p+1} \lg, W) \]
defined by equation \ref{differential}. 
Put $C^0=W$ and define $d^0(w)(\delta)=\nabla(\delta)(w)$. Let
$R_\nabla =
d^1\circ d^0$ be the \emph{curvature} of the connection $\nabla$.
Then from Lemma \ref{complex} it follows that we get a sequence  of
maps of $k$-vector spaces.
\end{definition} 

The Lie-Rinehart complex  is sometimes referred to as the \emph{Chevalley-Hochschild complex}.
We see from Lemma \ref{complex} that $C^\bullet (\lg ,W, \nabla )$ is a complex if and only if the curvature $R_\nabla$ 
is zero, hence if the curvature $R_\nabla$ is zero, we get well defined
cohomology spaces.

\begin{definition} \label{cohomology} Assume $\lg$ is a
  Lie-Rinehart algebra and $(W,\nabla)$ is a
a flat $\lg$-connection $\nabla$. We define the cohmology of $(W,\nabla )$ as follows:
 \[ \H^p(\lg,W,\nabla )=\H^p(C^\bullet (\lg,W, \nabla)), \]
where $C^\bullet(\lg, W, \nabla )$ is the Lie-Rinehart complex.
\end{definition} 
The maps $d^p$ from \ref{standardcomplex} are $k$-linear, hence the abelian groups $\H^p(\lg, W, \nabla)$ are 
$k$-vector spaces. 
Note furthermore that the cohomology $\H^*(\lg, A, \nabla)$ depends on the choice of connection
$\nabla :\lg \rightarrow \End_k(A)$.

If the ring $A$ is a smooth $k$-algebra of finite type, i.e the module of 
differentials $\Omega^1_{A}$ is locally free
of finite rank, it follows that the Lie-Rinehart complex 
is isomorphic to
the \emph{algebraic de Rham complex}, 
hence the Lie-Rinehart complex  generalizes simultaneously the algebraic 
de Rham complex and the Chevalley-Eilenberg complex.

\begin{proposition}\label{Deri} Let $A$ be a $k$-algebra and $\lg$ a
  Lie-Rinehart  algebra. Let furthermore $(W,\nabla_1)$ and $(W,\nabla_2)$ 
be $A$-modules with  $\lg$-connections. There exists an exterior-product
\[ C^*(\lg,W)\otimes_A C^*(\lg,W')\rightarrow C^*(\lg,W\otimes_A W') \]
with the following property:
\begin{equation}\label{derivation}
 d(xy )=d(x)y+(-1)^px d(y), 
\end{equation}
for all elements $x$ in 
$\Hom_A(\wedge^p\lg,W)$ and $y$ in $\Hom_A(\wedge^q\lg,W)$.
\end{proposition}
\begin{proof} Standard fact. 
\end{proof}

Recall some general definitions and standard facts on DG-algebras (for
a reference see \cite{weibel}).
A DG-algebra $B^*=\oplus_{p\geq 0}B^p$ is a graded associative algebra equipped
with a graded derivation $d$ of degree 1. If we do not require $d^2=0$
we say that $B^*$ is a \emph{quasi-differential graded algebra}. 
We can define the \emph{cohomology} $\H^*(B^*)$ of $B^*$ and it
follows that $\H^*(B^*)$ is a graded associative $k$-algebra. 
If $B^*$ is graded commutative, so is $\H^*(B^*)$. 
A graded left $B^*$-module $M^*=\oplus_{p\geq 0}M^p$ is a
\emph{differential graded module} if it is equipped with a graded
derivation of degree one with $d^2=0$. We say that $M^*$ is a
\emph{quasi-differential graded module} if we do not require $d^2=0$. 
It follows that $\H^*(M^*)$ is a graded left $\H^*(B^*)$-module. 
If we are given two DG-algebras $B^*$ and $E^*$ over a field $k$,
then a \emph{map of DG-algebras}, is just a map $\phi^*:B^*\rightarrow
E^*$ of graded associative algebras, commuting with the differentials.
One easily verifies that such a map $\phi^*$ induces
a map of graded associative algebras $\H(\phi^*):\H^*(B^*)\rightarrow 
\H^*(E^*)$. Also, given two DG-modules $M^*$ and $N^*$ on a DG-algebra
$B^*$, a \emph{map of DG-modules}, is a map $\psi^*:M^*\rightarrow
N^*$ commuting with the differentials. It is trivial to check
that such a map $\psi^*$ induces a map
$\H(\psi^*):\H^*(M^*)\rightarrow \H^*(N^*)$ of graded
$\H^*(B^*)$-modules.

\begin{proposition} Let $A$ be a $k$-algebra, $\mathfrak{g}$ a
Lie-Rinehart algebra. Let furthermore $(W,\nabla)$ be an $A$-module
with a $\lg$-connection. Then $C^*(\lg,A)$ is a DG-algebra and $C^*(\lg,W)$ is a
quasi-DG-module on $C^*(\lg,A)$. If $W$ is a $\lg$-module, then
$C^*(\lg,W)$ is a DG-module, hence $\H^*(\lg,A)$ is a graded  associative $k$-algebra
and $\H^*(\lg,W)$ is in a natural way a graded left module on $\H^*(\lg,A)$.
\end{proposition}
\begin{proof} This follows from the previous discussion and Proposition \ref{Deri}.
\end{proof}
\begin{proposition} Let $A$ be a $k$-algebra, and $\lg$ an
Lie-Rinehart algebra. Let furthermore $(W,\nabla )$ be a $\lg$-connection.
The connection  $\nabla $ induces a connection $ad\nabla $ on
$\End_A(W)$, 
hence $C^*(\lg,\End_A(W))$ becomes in a natural way a quasi-DG-algebra.
If $W$ is a $\lg$-module, $C^*(\lg,\End_A(W))$ is  DG-algebra.
\end{proposition}
\begin{proof} This follows from the previous discussion. 
\end{proof}

\section{A construction of the Chern character}

This section contains proofs of the following results: Let $k$ be a
field of characteristic zero, and let $A$ be a $k$-algebra. Let
furthermore $\mf{g}$ be an $(k,A)$-Lie-Rinehart algebra, and $(W,\nabla)$ be a
$\mf{g}$-connection wich is of finite presentation as an $A$-module. 
There exists a \emph{Chern character} $ch^\lg (W)$ in
$\H^*(\mf{g}|_U,\O_U)$ where $U$ is the open subset of $\Spec (A)$
where $W$ is locally free. We apply this to prove the existence of a ring homomorphism
\[ ch^\lg :\K_0(\lg)\rightarrow \H^*(\lg,A) \]
where $\K_0(\lg)$ is the \emph{Grothendieck ring} of locally free
$A$-modules with a $\lg$-connection.

 Recall briefly classical Chern-Weil theory: Let $A$ be a
$k$-algebra, where $k$ is a field of characteristic $0$, and let
$E$ be a locally free $A$-module. Any connection 
\[ \nabla : E \rightarrow \Omega^1_A\otimes E \] gives rise to a
connection
\[ ad\nabla :\End_A(E)\rightarrow \Omega^1_A\otimes \End_A(E)  ,\] 
and we get
\[ R_\nabla^k\in \Omega^{2k}_A\otimes \End_A(E). \]
Since $E$ is locally free there exists a trace map $tr:
\End_A(E)\rightarrow A$ and we get Chern-classes
\[ ch_k(E,\nabla)\in \H^{2k}_{\DR}(A). \]
This construction defines a group-homomorphism
\[ ch^A: \K_0(A)\rightarrow \H^*_{\DR}(A) .\]

\begin{theorem} \label{classicalchern} The map $ch^A:\K _0(A)\rightarrow \H^*_{\DR}(A)$ is a
  ring-homomorphism.
\end{theorem}
\begin{proof} See Theorem  8.1.7 in \cite{LODAY}.
\end{proof}
 
Notice the following: If $\nabla$ and $\nabla'$ are two
$\lg$-connections on an $A$-module $A$, where $\lg$ is an
Lie-Rinehart algebra,
then the difference $\nabla - \nabla'$ is an element of the module
$\Hom_A(\lg,\End_A(W))$. We express this by saying that the set
of $\lg$-connections on $W$ form a \emph{principal homogeneous space}
(or a \emph{torsor}) on $\Hom_A(\lg, \End_A(W))$. This means that
given a $\lg$-connection $\nabla$ on $W$, any other connection
$\nabla'$ can be obtained from $\nabla$ by adding an element
$\phi$ from $\Hom_A(\lg,\End_A(W))$, that is $\nabla'=\nabla + \phi$
for a unique $\phi$.

\begin{lemma}\label{commutative} Let $A$ be a $k$-algebra, $\lg$ a
  Lie-Rinehart algebra and $W$ a $\lg$-connection which is free as an $A$-module. The trace map
\[ tr^*:C^*(\lg, \End_A(W))\rightarrow C^*(\lg,A)\]
is a morphism of complexes.
\end{lemma}
\begin{proof} The only thing we have to
prove is that for all $p\geq 0$ we have commutative diagrams
\[ \diagram C^p(\lg, \End_A(W)) \rto^{d^p} \dto^{tr} &
C^{p+1}(\lg,\End_A(W)) \dto^{tr} \\
C^p(\lg,A) \rto^{d^p} & C^{p+1}(\lg,A) 
\enddiagram
\]
We may
assume that we have chosen a basis $\{e_i\}$ for $W$ as an $A$-module
and  we can write $W=\oplus_{i=1}^nAe_i$. Then in this basis we have
a connection $\nabla'_{\delta_i}(\sum a_ie_i)=\sum
\alpha(\delta_i)(a_i)e_i$,
and one verifies that $R_{\nabla'}=0$, hence the connection $\nabla '$
is integrable. The connection $\nabla$ which defines the $\lg$-structure
structure on $W$ can now be written in a unique way as $\nabla =\nabla
'+\phi$,
where $\phi$ is an element of $\Hom_A(\lg, \End_A(W))$, since 
$\lg$-connections are a torsor on $\Hom_A(\lg, \End_A(W))$. The
induced connection $ad\nabla$ on $\End_A(W)$ then becomes
\[ ad\nabla =[\nabla , - ]=[\nabla'+\phi,-]=[\nabla',-]+[\phi,-] .\]
The rest is straightforward calculation:
Let $\psi$ be an element of $C^p(\lg,\End_A(W))=\Hom_A(\wedge^p
\lg,\End_A(W))$. Put also $\omega =\delta_1 \wedge \cdots \wedge
\delta_{p+1}$, $\omega(i)=\delta_1\wedge \cdots \wedge
\hat{\delta_i}\wedge \cdots \wedge \delta_{p+1}$ for $1\leq i \leq
p+1$, and $\omega(i,j)=[\delta_i,\delta_j]\wedge \delta_1\wedge \cdots
\wedge \hat{\delta_i}\wedge \cdots \wedge \hat{\delta_j} \wedge \cdots
\wedge \delta_{p+1}$.

Then we have that
\[ tr((d^p\psi)(\omega)) \]
\[= tr\left(
  \sum_{i=1}^{p+1}(-1)^{i+1}ad\nabla_{\delta_i}\psi(\omega(i)) \right)
  + tr \left( \sum_{1\leq i < j \leq p+1}(-1)^{i+j}\psi(\omega(i,j))
\right)\]

\[=tr \left( \sum_{i=1}^{p+1}(-1)^{i+1}[\nabla'_{\delta_i},\psi(\omega(i))]
\right) + tr \left(
  \sum_{i=1}^{p+1}(-1)^{i+1}[\phi(\delta_i),\psi(\omega(i))] \right)
  \]
\[+tr \left( \sum_{1\leq i < j \leq p+1}(-1)^{i+j}\psi(\omega(i,j))
\right) \]
\[= tr\left( \sum_{i=1}^{p+1}(-1)^{i+1}(\alpha
  (\delta_i)(\psi(\omega(i)))_{kl} ) \right) 
 + tr\left( \sum_{i<j}(-1)^{i+j}\psi(\omega(i,j) ) \right) \] 
\[ =\sum_{i=1}^{p+1}(-1)^{i+1}tr(\alpha(\delta_i)(\psi(\omega(i))_{kl})
 ) + \sum_{i<j}(-1)^{i+j}tr(\psi(\omega(i,j))) \]
\[ =\sum_{i=1}^{p+1}(-1)^{i+1}\alpha(\delta_i)(tr(\psi(\omega(i)))) +
\sum_{i<j}(-1)^{i+j}(tr \circ \psi)(\omega(i,j)) \]
\[ =d^p(tr \circ \psi)(\omega)\]
and we see that $tr\circ d^p=d^p\circ tr$ and we have proved the assertion.
\end{proof}

\begin{corollary} \label{cohomologyclass}Assume $x^*$ is an element of $C^*(\lg,\End_A(W))$
  with the property that $d^*(x^*)=0$, then $tr^*(x^*)$ gives rise
to a cohomology-class $\overline{tr^*(x^*)}$ in $\H^*(\lg,A)$.
\end{corollary}
\begin{proof} We show that $d^*(tr^*(x^*))=0$: For all $p\geq 0$ 
we have commutative diagrams 
\[ \diagram C^p(\lg, \End_A(W)) \rto^{d^p} \dto^{tr} &
C^{p+1}(\lg,\End_A(W)) \dto^{tr} \\
C^p(\lg,A) \rto^{d^p} & C^{p+1}(\lg,A) 
\enddiagram
\]
by lemma \ref{commutative}. We see that
$d^p(tr^p(x^p)=tr(d^p(x^p))=tr(0)=0$,
hence we have that $d^*(tr^*(x^*))=0$ and we get a well-defined
cohomology-class $\overline{tr^*(x^*)}$ in $\H^*(\lg,A)$.
\end{proof}
Given a $\lg$-connection $W$, where $\lg$ is an Lie-Rinehart algebra,
one verifies that the curvature $R_\nabla$
is an element of $\Hom_A(\lg\wedge \lg,\End_A(W))=C^2(\lg,\End_A(W))$. 

\begin{lemma}(The Bianchi identiy)\label{bianchi} Let $A$ be a
$k$-algebra, $\lg$ a Lie-Rinehart algebra,
and $W$ a $\lg$-connection. Then $d^2(R_\nabla)=0$. 
\end{lemma}
\begin{proof} This is straightforward calculation:
Let $$\omega=\da\wedge \db\wedge \dc$$ be an element of
$\wedge^3\lg$.
Then we see that 
\[ d^2R_\nabla(\da\wedge\db\wedge \dc) \]
\[= ad\nabla_{\da}R_\nabla(\db\wedge \dc)-ad\nabla_{\db}
R_\nabla(\da\wedge \dc) +ad\nabla_{\dc}R_\nabla(\da\wedge
\db)\]
\[- R_\nabla([\da,\db]\wedge \dc)+R_\nabla([\da,\dc]\wedge
\db)-R_\nabla([\db,\da]\wedge \da) \]

\[= \nabla_{\da}R_\nabla(\db\wedge \dc)-R_\nabla(\db\wedge \dc)
\nabla_{\da} \]
\[- (\nabla_{\db}R_\nabla(\da\wedge \dc)-R_\nabla(\da\wedge
\dc)\nabla_{\db} ) +\nabla_{\dc}R_\nabla(\da\wedge \db)-
R_\nabla(\da\wedge \db)\nabla_{\dc} \]

\[- ([ \nabla_{[\da,\db]},\nabla_{\dc}
]-\nabla_{[[\da,\db],\dc]})
+[\nabla_{[\da,\dc]},\nabla_{\db}]-\nabla_{[[\da,\dc],\db]}\] 

\[-([\nabla_{[\db,\dc]},\nabla_{\da}]-\nabla_{[[\db,\dc],\da]})\]

\[=
[\nabla_{\da},[\nabla_{\db},\nabla_{\dc}]]+
[\nabla_{\db},[\nabla_{\dc},
\nabla_{\da}]]+[\nabla_{\dc},[\nabla_{\da},\nabla_{\db}]] \]
\[+ \nabla_{[[\da,\db],\dc]}+\nabla_{[[\db,\dc],\da]}+
\nabla_{[[\dc,\da],\db]} \]

\[+
[\nabla_{\db},\nabla_{[\da,\dc]}]-[\nabla_{\da},\nabla_{[\db,\dc]}]
\]
\[- [\nabla_{\dc},\nabla_{[\da,\db]}]\]

\[-
[\nabla_{[\da,\db]},\nabla_{\dc}]+[\nabla_{[\da,\dc]},\nabla_{\db}]
-[\nabla_{[\db,\dc]},\nabla_{\da}  ] =0\]
and we have proved the assertion.
\end{proof}

\begin{proposition}\label{dzero} Let $A$ be a $k$-algebra, $\lg$ a
  Lie-Rinehart algebra 
  and $(W,\nabla)$ be a $\lg$-connection. Let furthermore $R_\nabla$ be the
  curvature of $\nabla$. Then $d^{2n}(R_\nabla^n)=0$ for all $n\geq
  1$.
\end{proposition}
\begin{proof} We prove this by induction on $n$: By lemma
  \ref{bianchi} we see that the lemma is true for $n=1$. Assume it is
  true for $n=k$. We see that
\[ d(R_\nabla^{k+1})=d(R_\nabla^k\wedge R_\nabla)=d(R_\nabla^k)\wedge
  R_\nabla +(-1)^{2k}R_\nabla^k\wedge d(R_\nabla) \]
and $d(R_\nabla ^k)\wedge R_\nabla +(-1)^{2k}R_\nabla^k\wedge d(R_\nabla)$ 
is zero by the
induction hypothesis and lemma \ref{bianchi}, and we have proved the
assertion.
\end{proof}

Let in the following $A$ be a $k$-algebra, where $k$ is a field
of characteristic 0. Let $\lg$ be a Lie-Rinehart algebra and
$(W,\nabla )$ a $\lg$-connection, where $W$ is an $A$-module of finite
presentation. Let $exp(R_\nabla)$ be defined as $\sum_{n\geq
  0}\frac{1}{n!}R_\nabla^n$. Consider the open set $U\subseteq \Spec A$ 
where $W$ is locally free, which exists since $W$ is of finite
presentation. By lemma \ref{commutative} we have trace maps
\[ tr^*:C^p(\lg_\mf{p},\End_{A_\mf{p}}(W_\mf{p}))\rightarrow
C^*(\lg_\mf{p},A_\mf{p}) \]
defined for all $\mf{p}$ in $U$, since $W|_U$ is locally free, and
these maps glue to give a map of sheaves of complexes
\[ tr^*: C^*(\lg|_U,\End_{\O_U}(W|_U))\rightarrow C^*(\lg|_U,\O_U) .\]
We have that $R_\nabla|_U$ is an element of $C^2(\lg|_U,
\End_{\O_U}(W|_U))$
and we obtain an element $exp(R_\nabla|_U)$ in
$C^*(\lg|_U,\End_{\O_U}(W|_U))$. By lemma \ref{dzero} we see  that
the element
$d^*(exp(R_\nabla|_U))$ equals zero, since it vanishes when we localize at all
prime-ideals $\mf{p}$ in $U$. Consider the element
$x^*=tr^*(exp(R_\nabla|_U))$, which lives in $C^*(\lg|_U,\O_U)$.

\begin{theorem} \label{theoremone} The following holde: $d^*(x^*)=0$. Hence $\overline{x^*}$
  defines a cohomology-class in $\H^*(\lg|_U,\O_U)$. 
\end{theorem}
\begin{proof} It follows  from corollary \ref{cohomologyclass} that $d^*(x^*)=0$, since
  we have already seen that $d^*(exp(R_\nabla))=0$, hence we get a
cohomology class as claimed. 
\end{proof}

\begin{definition} Let $A$ be a $k$-algebra where $k$ is a field of
  characteristic 0 and let $\lg$ be an Lie-Rinehart algebra. Let
  furthermore $W$ be a $\lg$-connection, where $W$ is an
  $A$-module  of finite presentation. We let the element 
  $ch^\lg (W, \nabla)=\overline{x^*}$ in $\H^*(\lg|_U,\O_U)$ from theorem \ref{theoremone}
  be the \emph{Chern character} of the $\lg$-connection $(W,\nabla)$.
\end{definition}

By theorem \ref{theoremone} the class $ch^\lg (W,\nabla )$ in
$\H^*(\lg|_U,\O_U)$ is an invariant of the pair $(W,\nabla)$.
Given any $k$-algebra $A$, where $k$ is a field of characteristic 0,
and $\lg$ an Lie-Rinehart algebra, we consider $\K_0(\lg)$, the
\emph{Grothendieck ring} of $\lg$. This is defined as follows:
$\K_0(\lg)$ is the free abelian group on the symbols $[W,\nabla]$ module a subgroup
$D$ wich we will define below.
Here $(W,\nabla)$ is a  $\lg$-connection which is a locally free $A$-module of finite
rank. The 
symbol $[W,\nabla]$ denotes the isomorphism-class of the pair
$(W,\nabla)$. The subgroup $D$ is the group 
generated by the relations
\[ [W\oplus W',\nabla \oplus\nabla']-[W,\nabla]-[W',\nabla'] .\]
That is: $\K_0(\lg)=\oplus \mathbf{Z}[W,\nabla]/D$.
(We obviously have that the direct sum of two $\lg$-connections is again a
$\lg$-connection.) 
Given two $\lg$-connections $(W,\nabla)$ and $(W',\nabla')$, there
exists a natural connection $\nabla \otimes \nabla'=\nabla \otimes 1 +
1\otimes \nabla'$ on $W\otimes _A W'$, hence $W\otimes _A W'$ is in a natural
way a $\lg$-connection. Define a map
\[\otimes: \oplus \mathbf{Z}[W,\nabla ] \times \oplus
\mathbf{Z}[W,\nabla ] \rightarrow \K_0(\lg) \]
by the following 
\[  \otimes(\sum_i n_i[W_i,\nabla _i],\sum_j m_j[V_j,\nabla' _j])=
\sum_{i,j}n_im_j[W_i\otimes_A V_j,\nabla _i\otimes \nabla '_j]. \]

\begin{lemma} The map $\otimes$ defines a $\mathbf{Z}$-bilinear product
\[ \K_0(\lg)\otimes_{\mathbf{Z} }\K_0(\lg)\rightarrow \K_0(\lg) \]
making $\K_0(\lg)$ into a commutative $\mathbf{Z}$-algebra.
\end{lemma}
\begin{proof} This is straightforward. 
\end{proof}

\begin{lemma} \label{connectionprops} 
Let $(W,\nabla)$ and $(W',\nabla')$ be two $\lg$-connections.
The the following holds:
\begin{align}
& R_{\nabla \oplus \nabla'}=R_\nabla \oplus R_{\nabla'} \label{r1}\\
& R_{\nabla \otimes \nabla'}=R_\nabla \otimes 1+1\otimes R_{\nabla '} \label{r2} \\
&R_\nabla \otimes 1 \wedge 1 \otimes R_{\nabla'}=1\otimes R_{\nabla'}\wedge
R_\nabla \otimes 1 \label{r3} \\
&(R_{\nabla \oplus \nabla'})^n=R_\nabla^n\oplus R_{\nabla'}^n \label{r4}
\end{align}
\end{lemma}
\begin{proof} 
We first prove equation \ref{r1}:

\[ R_{\nabla \oplus \nabla '}(\delta \wedge \eta)= [\ddd,\dde]-\dd
_{[\delta,\eta]} .\]
It follows that if we pick $(w,w')$ in $W\oplus W'$, we get
\[ R_{\dd }(\delta \wedge \eta)(w,w') \]
\[  =[\ddd ,\dde ](w,w')-\dd_{[\delta,\eta]}(w,w') \]
\[ =\dde \circ \ddd(w,w')-\ddd \circ \dde(w,w')-
   \dd_{[\delta,\eta]}(w,w') \]
\[= \dde (\nd (w),\nd ' (w') )-\ddd (\ne (w), \ne '(w') ) 
-(\nabla _{[\delta ,\eta ]}(w),\nabla _{[\delta ,\eta ]}'(w') ) \]
\[=(\ne \nd (w),\ne ' \nd '(w') )-(\nd \ne (w),\nd ' \ne '(w') ) 
  -(\nabla _{[\delta ,\eta ]}(w),\nabla _{[\delta ,\eta ]}'(w') ) \]
 \[  =(R_\nabla (\delta \wedge \eta )(w),R_{\nabla '}(\delta \wedge \eta )(w') ) \]
 \[=R_\nabla \oplus R_{\nabla '}(w,w') \]

and equation \ref{r1} follows. We prove equation \ref{r2}: Let 
$w\otimes w'$ be an element of $W\otimes_A W'$, and let $\dtd =
\nabla \otimes 1 + 1\otimes \nabla'$ be the $\lg$-connection on 
$W\otimes_A W'$. We get
\[ R_{\dtd } (\delta \wedge \eta)(w\otimes w') \]
\[=[\dtdd , \dtde ](w\otimes w')-\dtd _{[\delta , \eta ]}(w\otimes w') \]
\[=\dtde \circ \dtdd (w\otimes w')-\dtdd \circ \dtde (w\otimes w')
-\dtd _{[\delta , \eta ]}(w\otimes w') \]
\[=\dtde (\nd (w)\otimes w' +w\otimes \nd '(w') )-\dtdd (\ne (w)\otimes w'
+w\otimes \ne '(w') \]
\[ -(\nabla _{[\delta ,\eta ]}(w)\otimes w'+w\otimes \nabla _{[\delta ,\eta]}(w') ) \]
\[ =\ne \nd (w)\otimes w' + \nd (w)\otimes \ne '(w')+\ne (w)\otimes \nd '(w')
+w\otimes \ne '\nd '(w') \]
\[-(\nd \ne (w)\otimes w'+\ne (w)\otimes \nd '(w') +\nd (w)\otimes \ne '(w')
+w\otimes \nd ' \ne ' (w') )\]
\[ -\nabla_{[\delta ,\eta ]}(w)\otimes w' - w\otimes \nabla_{[\delta , \eta ]}
(w') \]
\[=[\nd , \ne ](w)\otimes w' + w\otimes [\nd ' , \ne '](w')
-\nabla _{[\delta , \eta ]}(w)\otimes w' -w\otimes \nabla'_{[\delta ,\eta ]}
(w') \]
\[=R_\nabla (\delta \wedge \eta)(w)\otimes w'+w\otimes R_{\nabla '}(\delta \wedge \eta )(w') \]
and equation  \ref{r2} follows. We prove equation \ref{r3}: Let
$\omega $ be an element of $\wedge ^4 \lg$. We get
\[ R_\nabla \otimes 1\wedge 1\otimes R_{\nabla '}(w)\]
\[=\sum_{(2,2)}sgn(\sigma )(R_\nabla \otimes 1, 1\otimes R_{\nabla '})
\sigma (\omega) \]
\[=\sum_{(2,2)}sgn( \sigma )R_\nabla (\sigma (\omega))\otimes 1 \circ
1\otimes R_{\nabla '}(\sigma (\omega )) \]
\[ =\sum_{(2,2)}sgn(\sigma )1\otimes R_{\nabla '}(\sigma (\omega))\circ
R_\nabla (\sigma (\omega))\otimes 1 \]
\[=\sum_{(2,2)}sgn(\sigma )(1\otimes R_{\nabla '},R_\nabla \otimes 1)
\sigma (\omega) \]
\[=1\otimes R_{\nabla '}\wedge R_{\nabla }\otimes 1 (\omega) \]
and equation \ref{r3} follows. Finally we prove equation \ref{r4} by induction
on $n$. For $n$=2 we get the following:
Let $\omega =\delta_1\wedge \cdots \wedge \delta_4$, and for any 
$(2,2)$-shuffle $\sigma$ put $\sigma (\omega )^1=\delta_{\sigma (1)}\wedge
\delta_{\sigma (2)}$ and $\sigma(\omega )^2=\delta_{\sigma (3)}\wedge
\delta_{\sigma (4)}$. We get
\[ (R_\nabla \oplus R_{\nabla '})^2(\omega) \]
\[=\sum_{(2,2)}sgn(\sigma )(R_\nabla \oplus R_{\nabla '},R_\nabla \oplus R_{\nabla '})\sigma (\omega) \]
\[=\sum_{(2,2)} sgn(\sigma )R_\nabla \oplus R_{\nabla '}(\sigma (\omega )^1)
\circ R_\nabla \oplus R_{\nabla '}(\sigma (\omega )^2) \]
\[=\sum_{(2,2)} sgn(\sigma )R_\nabla(\sigma (\omega)^1)\oplus
R_{\nabla '}(\sigma (\omega )^1)\circ R_\nabla (\sigma (\omega )^2)\oplus
R_{\nabla '}(\sigma (\omega )^2) \]
\[=\sum_{(2,2)} sgn(\sigma ) R_\nabla (\sigma (\omega)^1)R_\nabla (\sigma 
(\omega )^2)\oplus R_{\nabla '}(\sigma (\omega )^1)R_{\nabla '}(\sigma 
(\omega )^2) \]
\[=(\sum_{(2,2)}sgn(\sigma )(R_\nabla ,R_\nabla )\sigma (\omega ))\oplus
(\sum_{(2,2)}sgn(\sigma)(R_{\nabla '},R_{\nabla '})\sigma (\omega )) \]
\[=R_\nabla^2(\omega )\oplus R_{\nabla '}^2(\omega ) \]
and we have proved equation \ref{r4} for $n=2$.
Assume the equation is true for $n=k$. Put $n=k+1$, and let
$\omega =\delta_1 \wedge \cdots \wedge \delta_{2k+2}$. Put also
for any $(2k,2)$-shuffle $\sigma$, $\sigma (\omega )^1=\delta_{\sigma (1)}
\wedge \cdots \wedge \delta_{\sigma (2k)}$ and $\sigma (\omega )^2=
\delta_{\sigma (2k+1)}\wedge \delta_{\sigma (2k+2)}$. We get
\[ R_{\dd}^kR_{\dd }(\omega) \]
\[=\sum_{(2k,2)}sgn(\sigma )(R_{\dd}^k,R_{\dd })\sigma (\omega ). \]
By the induction hypothesis we get
\[=\sum_{(2k,2)}sgn(\sigma )(R_\nabla ^k\oplus R_{\nabla '}^k,R_\nabla 
\oplus R_{\nabla '})\sigma (\omega ) \]
\[=\sum_{(2k,2)}sgn(\sigma ) R_{\nabla }^k\oplus R_{\nabla '}^k(\sigma 
(\omega )^1)\circ R_\nabla \oplus R_{\nabla '}(\sigma (\omega )^2) \]
\[=\sum_{(2k,2)}sgn(\sigma )R_{\nabla }^k(\sigma (\omega )^1)\oplus
R_{\nabla '}^k(\sigma (\omega )^1)\circ R_\nabla (\sigma (\omega )^2)\oplus
R_{\nabla '}(\sigma (\omega )^2) \]
\[=\sum_{(2k,2)}sgn(\sigma )R_\nabla ^k(\sigma (\omega )^1)R_\nabla
(\sigma (\omega )^2)\oplus R_{\nabla '}^k(\sigma (\omega )1)\circ 
R_{\nabla '}(\sigma (\omega )^2) \]
\[=(\sum_{(2k,2)}sgn(\sigma )(R_\nabla ^k,R_\nabla )\sigma (\omega ))\oplus
(\sum_{(2k,2)} sgn(\sigma )(R_{\nabla '}^k,R_{\nabla '})\sigma (\omega )) \]
\[=(R_\nabla ^{k+1}\oplus R_{\nabla '}^{k+1})(\omega ) \]
and equation \ref{r4} follows, and we have proved the lemma.
\end{proof}

\begin{lemma} \label{tracetensor}
Let $W$ and $W'$ be two free $A$-modules, and let
$\phi$ in $\End_A(W)$ and $\psi$ in $\End_A(W')$ be two endomorphisms.
Then the following holds 
\[ tr(\phi \otimes \psi)=tr(\phi)tr(\psi). \]
\end{lemma}
\begin{proof} Let $W=\oplus_{i=1}^n Ae_i$ and $W'=\oplus_{j=1}^mAf_j$ 
be two direct-sum decompositions of $W$ and $W'$. Put also
$\phi =(a_{ij})$ and $\psi=(b_{ij})$ where $a_{ij}$ and $b_{ij}$ are
elements of $A$. One verifies that for instance $tr(\phi)=\sum_i e_i\phi e_i$.
We get
\[ tr(\phi \otimes \psi)=\sum_{i,j}e_i\otimes f_j(\phi\otimes \psi)e_i\otimes f_j .\]
It is trivial to check that $e_k\otimes f_l(\phi\otimes \psi )e_m\otimes f_k$
equals $a_{km}b_{ln}$, hence we get
\[ \sum_{i,j}a_{ii}b_{jj}=(\sum_i a_{ii})(\sum_j b_{jj} )=(tr\phi)(tr\psi) \]
and the lemma follows.
\end{proof}

\begin{lemma}\label{tracewedge}
Let $(W,\nabla)$ and $(W',\nabla')$ be two  locally
free $\lg$-connections, then
\[ tr(R_\nabla^n\otimes 1 \wedge 1\otimes R_{\nabla'}^m )=(tr(R_\nabla^n))
\wedge (tr(R_{\nabla'}^m )) \]
\end{lemma}
\begin{proof} 
Let
$\omega =\delta_1 \wedge \cdots \wedge \delta_{2(n+m)}$, and put
for any $(2n,2m)$ shuffle $\sigma $, $\sigma (\omega )^1=
\delta_{\sigma (1)}\wedge \cdots \wedge \delta_{\sigma (2n) }$ and
$\sigma (\omega )^2=\delta_{\sigma (2n+1) }\wedge \cdots \wedge
\delta_{\sigma (2(n+m) )}$. We see that
\[ R_\nabla^n\otimes 1\wedge 1\otimes R_{\nabla'}^m(\omega) \]
\[= \sum_{(2n,2m)}sgn(\sigma )(R_\nabla^n\otimes 1,1\otimes R_{\nabla '}^m)
\sigma (\omega ) \]
\[ =\sum_{(2n,2m)}sgn(\sigma )R_\nabla^n(\sigma (\omega )^1)\otimes 1 \circ
1\otimes R_{\nabla '}^m(\sigma (\omega )^2) \]
\[=\sum_{(2n,2m)}sgn(\sigma ) R_\nabla^n(\sigma (\omega )^1)\otimes
R_{\nabla '}^m(\sigma (\omega )^2). \]
By lemma \ref{tracetensor} we get
\[ tr(R_\nabla^n\otimes 1\wedge 1\otimes R_{\nabla'}^m(\omega)) \]
\[=tr(\sum_{(2n,2m)}sgn(\sigma ) R_\nabla^n(\sigma (\omega )^1)\otimes
R_{\nabla '}^m(\sigma (\omega )^2) ) \]
\[=\sum_{(2n,2m)}sgn(\sigma )(tr\circ R_\nabla ^n)(\sigma (\omega )^1)(tr 
\circ \R_{\nabla '}^m)(\sigma (\omega )^2) \]
\[=\sum_{(2n,2m)}sgn(\sigma )(tr\circ R_\nabla^n,tr\circ R_{\nabla '}^m)
\sigma (\omega) =(tr\circ R_\nabla^n)\wedge (tr\circ R_{\nabla '}^m)(\omega)
\]
and we have proved the assertion.
\end{proof}

We can now prove the existence of the Chern character.

\begin{theorem} \label{theoremtwo} There exists a ring homomorphism
\[ ch^\lg :\K_0(\lg)\rightarrow \H^*(\lg,A) \]
from the Grothendieck ring  $\K_0(\lg)$ 
to the cohomology ring $\H^*(\lg,A)$.
\end{theorem}
\begin{proof} For every locally free $\lg$-connection $W$ of finite rank

we obtain by Theorem \ref{theoremone} a cohomology class

$ch^\lg (W)$ in $\H^*(\lg,A)$. 
Define a map  $\phi: \oplus \mathbf{Z}[W,\nabla]\rightarrow 
\H^*(\lg,A)$ by the formula
\[ \phi(\sum_i n_i[W_i,\nabla _i]=\sum_i n_ich(W_i, \nabla _i) .\]
We want to show that the map $\phi$ gives rise to  a well-defined map
\[ ch^\lg :\K_0(\lg)\rightarrow \H^*(\lg, A) .\]
Let $[W\oplus W',\nabla\oplus \nabla']-[W,\nabla ]-[W',\nabla']$ be a generator of the group
$D$, where $\K_0(\lg)=\oplus \mathbf{Z}[W,\nabla]/D$. We get
\[ ch^\lg ([W\oplus W',\nabla \oplus \nabla']-[W,\nabla
]-[W',\nabla'])= \]
\[ ch^\lg (W\oplus W',\nabla \oplus \nabla')-ch^\lg(W,\nabla )-ch^\lg(W',\nabla') \]
\[=\sum_{n\geq 0}\frac{1}{n!}tr(R_{\nabla \oplus \nabla'})^n
-\sum_{k\geq 0}\frac{1}{k!}trR_\nabla^k-\sum_{l\geq 0}\frac{1}{l!}trR
_{\nabla '}^l .\]
By lemma \ref{connectionprops}, equation \ref{r1} and \ref{r4}
we get
\[ \sum_{n\geq 0}\frac{1}{n!}tr(R_\nabla^n\oplus R_{\nabla '}^n)
-\sum_{k\geq 0}\frac{1}{k!}tr(R_\nabla^k)-\sum_{l\geq 0}\frac{1}{l!}
tr(R_{\nabla '}^l) \]
\[= \sum_{n\geq 0}\frac{1}{n!}(trR_{\nabla }^n+trR_{\nabla '}^n)
-\sum_{k\geq 0}\frac{1}{k!}trR_\nabla^k-\sum_{l\geq 0}trR_{\nabla '}^l =0 \]
hence $\phi$ gives rise to a map $ch^\lg :\K_0(\lg)\rightarrow \H^*(\lg,A)$, and
obviously $ch^\lg$ is a group-homomorphism. We show that $ch^\lg $ is a 
ring homomorphism: Put for any $\lg$-connection $(W,\nabla)$, $ch_n(W,\nabla)=
\frac{1}{n!}trR_\nabla^n$. We have that $ch^\lg (W,\nabla)=\sum_{n\geq 0}ch_n(W,\nabla)$. 
Since $C^*(\lg, \End_A(W\otimes_A W'))$ is an associative
$A$-algebra and by lemma \ref{connectionprops}, equation \ref{r3}
we have that $R_\nabla \otimes 1 \wedge 1\otimes R_{\nabla '}=
1\otimes R_{\nabla '}\wedge R_\nabla \otimes 1$, we can apply the
binomial-theorem. We get
\[ ch_n(W\otimes W',\nabla \otimes \nabla ')=\frac{1}{n!}(R_{\nabla \otimes \nabla '})^n \]
and by lemma \ref{connectionprops}, equation \ref{r2} we get
\[ \frac{1}{n!}tr(R_\nabla \otimes 1 +1\otimes R_{\nabla '})^n 
= \sum_{i+j=n}\frac{1}{i!j!}tr(R_\nabla\otimes 1)^i(1\otimes R_{\nabla '})^j 
.\] 
By lemma \ref{tracewedge} we get
\[ \sum_{i+j=n}\frac{1}{i!j!}(trR_\nabla)^i\wedge (trR_{\nabla '})^j =
\sum_{i+j=n}(\frac{1}{i!}trR_\nabla^i)\wedge(\frac{1}{j!}trR_{\nabla '}^j) \]
\[ =\sum_{i+j=n}ch_i(W,\nabla )ch_j(W',\nabla' ) .\]
The following holds
\[ ch^\lg (W\otimes W',\nabla \otimes \nabla')=\sum_{n\geq 0}ch_n(W\otimes
W',\nabla \otimes \nabla' )\]
\[=\sum_{n\geq 0}\left( \sum_{i+j=n}ch_i(W,\nabla )ch_j(W',\nabla' ) \right)=
\left(\sum_{k\geq 0}ch_k(W,\nabla ) \right)\left(\sum_{l\geq
    0}ch_l(W',\nabla' ) \right) \]
\[=ch^\lg (W,\nabla )ch^\lg (W',\nabla' ), \]
and the theorem follows.
\end{proof}

\section{On independence of choice of connection}

In this section we prove 
the fact that the Chern character $ch^\lg (W,\nabla)$ of
an $A$-module with a $\lg$-connection from Theorem \ref{theoremone} is independent with respect to choice of
connection $\nabla$. 
Let in the following $A$ be a $k$-algebra where $k$ is a field of
characteristic zero. Let furthermore $\lg$ be a Lie-Rinehart algebra with
anchor map $\alpha:\lg \rightarrow \Der_k(A)$.
We first prove a series of technical lemmas: 

\begin{lemma} We get in a natural way a map $\alpha\otimes
  1:\lg[t]\rightarrow \Der_k(A[t])$, 
making $\lg[t]$ into an $(k,A[t])$-Lie-Rinehart algebra.
\end{lemma}
\begin{proof} Define a $k$-Lie algebra structure on $\lg [t]$ as
  follows: $[\sum_i \delta _i\otimes f_j,\sum_j \eta_j\otimes
  g_j]=\sum_{i,j}[\delta_i,\eta_j]\otimes f_ig_j$. Define furthermore
  a map $\alpha\otimes 1:\lg [t]\rightarrow \Der_k(A[t])$ by
$\alpha \otimes 1(\delta\otimes f)(a\otimes
  g)=\alpha(\delta)(a)\otimes fg$, then it is straightforward to
  check that $\lg [t]$ is a $(k,A[t])$-Lie-Rinehart algebra.
\end{proof}

\begin{lemma} Let $W$ be an $A$-module with a $\lg$-connection
  $\nabla$. There exists a $\lg[t]$-connection $\nabla \otimes 1$ on the $A[t]$-module $W[t]$.
\end{lemma}
\begin{proof} Define the following map: $\nabla \otimes 1:\lg
  [t]\rightarrow \End_k(W[t])$, by letting $\nabla \otimes 1(\delta
  \otimes f)(w\otimes g)=\nabla (\delta)(w)\otimes fg$. Then it is
  straightforward to check that $\nabla \otimes 1$ is a $\lg [t]$-connection.
\end{proof}

\begin{lemma} Let $\nabla_0$ and $\nabla_1$ be $\lg$-connections on
  $W$, then $\nabla=\nabla_1\otimes t+\nabla_0\otimes (1-t)$ is a
  $\lg[t]$-connection on $W[t]$.
\end{lemma}
\begin{proof} This is straightforward.
\end{proof}

\begin{lemma} Let $\nabla$ be a $\lg$-connection on an $A$-module
  $W$. Let $\nabla \otimes 1$ be the induced $\lg[t]$-connection on
  $W[t]$. Then the curvature $R_{\nabla \otimes 1}$ defines a natural
  map
\[ R_{\nabla \otimes 1}:\wedge^2 \lg[t]\rightarrow \End_A(W)[t] .\]
\end{lemma}
\begin{proof} Define $R_{\nabla \otimes 1}(\delta \otimes f \wedge
  \eta\otimes g)=R_\nabla (\delta \wedge \eta)\otimes fg$, then the
  lemma follows. 
\end{proof}

\begin{lemma} Let $\nabla$ be a $\lg$-connection on the $A$-module
  $W$, and consider the induced connection $\nabla\otimes 1$ on
  $W[t]$. There exists
a map $p^i_*:C^p(\lg[t],W[t])\rightarrow C^p(\lg, W)$ making
  commutative diagrams
\[
\diagram C^p(\lg[t],W[t]) \rto^d \dto^{p^i_*} & C^{p+1}(\lg[t],W[t]) \dto^{p^i_*}
\\
         C^p(\lg,W) \rto^d & C^{p+1}(\lg,W) 
\enddiagram \]
for all $p$.
\end{lemma}
\begin{proof} Define the maps 
\[ p^i_*: C^p(\lg [t],W[t])\rightarrow C^{p}(\lg [t],W[t]) \]
as follows: There exists an obvious map $q:\wedge ^p\lg \rightarrow
\wedge ^p \lg [t]$ defined by mapping $\delta_1\wedge \cdots \wedge \delta_p$
to $\delta_1\otimes 1\wedge \cdots \wedge \delta_p\otimes 1$. There
exists a map $p^i: W[t]\rightarrow W$ defined by letting $p^i(t)=i$
for $i=0,1$. Put now for any $A$-linear map $\phi:\wedge ^p \lg
[t]\rightarrow W[t]$, $p^i_*(\phi)=p^i\circ \phi \circ q$. We show
that we get commutative diagrams as claimed: Consider first
$p^i_*(d\phi)(\delta_1\wedge \cdots \wedge \delta_{p+1})=$
\[ p^id(\phi)(\delta_1\otimes 1\wedge \cdots \wedge
\delta_{p+1}\otimes 1)= \]
\[ p^i( \sum_{k=1}^{p+1}(-1)^{k+1}\nabla (\delta_k)\otimes 1
  \phi(\cdots \wedge \hat{\delta_k \otimes 1} \wedge \cdots) + \]
\[ \sum_{k<l}(-1)^{k+l}\phi([\delta_k\otimes 1, \delta_l \otimes
1]\wedge \cdots \hat{\delta_k\otimes 1} \cdots \hat{\delta_l \otimes
  1} \cdots ) )= \]
\begin{equation} \label{c1}
 \sum_{k=1}^{p+1}p^i\nabla (\delta_k)\otimes 1 \phi(\cdots
  \hat{\delta_k \otimes 1} \cdots ) + 
\end{equation}
\[\sum_{k<l}(-1)^{k+l}p^i\phi([\delta_k,\delta_l]\otimes 1\wedge
\cdots \hat{\delta_k \otimes 1}\cdots \hat{\delta_l \otimes 1} \cdots)
.\]
Consider $d(p^i_*\phi)(\delta_1\wedge \cdots \wedge \delta_{p+!}) =$
\[ \sum_{k=1}^{p+1}(-1)^{k+1}\nabla (\delta_k)p^i_*\phi(\cdots
\hat{\delta_k} \cdots ) + \]
\[ \sum_{k<l}(-1)^{k+l}p^i_*\phi([\delta_k,\delta_l]\cdots
\hat{\delta_k} \cdots \hat{\delta_l} \cdots ) =\]
\begin{equation} \label{c2}
\sum_{k=1}^{p+1}(-1)^{k+1}\nabla (\delta_k)p^i\phi(\cdots \hat{
  \delta_1 \otimes 1 } \cdots ) + 
\end{equation}

\[ \sum_{k<l}(-1)^{k+l} p^i\phi([\delta_k,\delta_l]\otimes 1  \cdots
  \hat{\delta_k\otimes 1 } \cdots \hat{\delta_l \otimes 1 } \cdots )
  .\]
One checks that $\nabla (\delta_k)p^i=p^i \nabla(\delta_k)\otimes 1$
hence equation \ref{c1} equals equation \ref{c2}, and the claim follows.
\end{proof}

\begin{lemma} Given two $\lg$-connections $\nabla_0,\nabla_1$ on $W$,
  and let $\nabla=\nabla_1\otimes t+\nabla_0\otimes (1-t)$ be the
  induced connection on $W[t]$. Then the curvature $R_{\nabla}$ 
is an element of $C^2(\lg[t],\End_A(W)[t])$, and it follows that 
$p^i_*(R_{\nabla})=R_{\nabla_i}$ for $i=0$ and $1$.
\end{lemma}
\begin{proof} This is straighforward. \end{proof} 

\begin{lemma} Consider the map 
\[ p^i_*:C^p(\lg[t],\End_A(W)[t])\rightarrow C^p(\lg,\End_A(W)) .\]
Let $\phi$ and $\psi$ be elements of $C^p(\lg[t],\End_A(W)[t])$ and
$C^q(\lg[t],\End_A(W)[t])$ respectively. The following holds: 
\[p^i_*(\phi \wedge \psi)=p^i_*(\phi)\wedge p^i_*(\psi).\]
In particular it follows that $p^i_*(R_\nabla^k)=(p^i_*R_\nabla)^k$. 
\end{lemma}
\begin{proof} This is straighforward. \end{proof} 

\begin{lemma} There exists for all $p$ commutative diagrams
\[ 
\diagram C^p(\lg[t],\End_A(W)[t])\rto^{tr\otimes 1}\dto^{p^i_*} &
C^p(\lg[t],A[t]) \dto^{p^i_*} \\
C^p(\lg, \End_A(W))\rto^{tr} & C^p(\lg,A)  
\enddiagram \]
in particular we get
$p^i_*(tr(R_\nabla^k))=tr(p^i_*R_\nabla^k)$.
\end{lemma}
\begin{proof} Let $\phi :\wedge^p \lg \rightarrow \End_A(W)[t]$ be an
  $A$-linear map. Since $W$ is locally free, we have a trace map
$tr:\End_A(W)\rightarrow A$, and we get a trace-map $tr\otimes 1:
\End_A(W)[t]\rightarrow A[t]$, and we get $tr\otimes 1 \circ \phi$ in 
$C^p(\lg [t],A[t])$. We see that $p^i_*(tr\otimes 1 \circ
\phi)(\delta_1\wedge 
\cdots \wedge \delta_p)=$
\begin{equation} \label{co1}
 p^i_*\circ tr\otimes 1 \circ \phi(\delta_1\otimes 1
  \wedge \cdots \wedge \delta_p\otimes 1 ) 
\end{equation}
We also see that $tr(p^i_*(\phi))(\delta_1\wedge \cdots \wedge
\delta_p)=$
\begin{equation} \label{co2}
tr\circ p^i \circ \phi(\delta_1\otimes 1 \wedge
  \cdots \wedge \delta_p\otimes 1)
\end{equation}
and since $p^i_*\circ tr\otimes 1 =tr \circ p^i_*$
we see that equation \ref{co1} equals equation \ref{co2}, and we have
proved the assertion.
\end{proof}

\begin{lemma} The maps $p^i_*:C^p(\lg [t], \End_A(W)[t])\rightarrow
  C^p(\lg, \End_A(W))$ satisfy $p^i_*(\phi \wedge
  \psi)=p^i_*(\phi)\wedge p^i_*(\psi)$. In particular we get
 $p^i_*(R_\nabla^k)=(p^i_*R_\nabla)^k$.
\end{lemma}
\begin{proof}  $p^i_*(\phi \wedge \psi)(\delta_1 \wedge \cdots \wedge \delta_{p+q})=$
\[ p^i(\phi \wedge \psi(\delta_1\otimes 1 \wedge \cdots \wedge
  \delta_{p+q}\otimes 1) =\]
\[p_i \sum_{(p,q)} sgn(\sigma )\phi(\delta_{\sigma (1)}\otimes 1
\cdots \delta_{\sigma (p)}\otimes 1 )\psi(\delta_{\sigma (p+1)}\otimes
1 \cdots \delta_{\sigma(p+q)}\otimes 1 )= \]
\[ \sum_{(p,q)} sgn( \sigma)p^i_*\phi(\delta_{\sigma(1)}\cdots
\delta_{\sigma(p)} )p^i_*(\psi)(\delta_{\sigma(p+1)}\cdots
\delta_{\sigma(p+q)}) =\]
\[ p^i_*(\phi)\wedge p^i_*(\psi)(\delta_1 \wedge \cdots \wedge
\delta_{p+q}) \]
and the lemma follows.
\end{proof}

We are now in position to prove the main theorem of this section. 

\begin{theorem}\label{independence}  Let $A$ be any $k$-algebra where $k$ is any field,
and let  $\lg$ be a Lie-Rinehart algebra. Let $W$ be a locally free
  $A$-module with a $\lg$-connection $\nabla$. The class
$ch_n(W,\nabla)$ in $\H^{2n}(\lg,A)$ is independent with respect to choice
of connection.
\end{theorem}
\begin{proof} Consider the complex $C^*(\lg [t], A[t])$:

\[ \cdots \rightarrow C^{p-1}(\lg[t], A[t]) \rightarrow C^p(\lg[t],A[t]) \rightarrow C^{p+1}(\lg[t],A[t]) \rightarrow \cdots \]
By functoriality we get:
\[ C^p(\lg[t],A[t])=\Hom_A(\wedge^p(\lg\otimes_A A[t]),A[t])=\Hom_A((\wedge^p \lg)\otimes _A A[t],A[t])= \]
\[ \Hom_A(\wedge^p \lg, A)\otimes _A A[t]=\Hom_A(\wedge^p \lg ,A)\otimes _k k[t] .\]
It follows that we get an isomorphism
at the level of cohomology-groups 
\[ \H^i(\lg [t],A[t])\cong \H^i(\lg, A) .\]
We get induced maps on cohomology groups
\[ p^i_*:\H^{2k}(\lg [t], A[t])\rightarrow \H^{2k}(\lg ,A) \]
with the property that 
\[ p^i_*(\overline{tr(R_\nabla)} )=\overline{tr(R_{\nabla_i})}.\]
It follows that 
\[ \overline{tr(R_{\nabla_0})}=\overline{tr(R_{\nabla_1})} ,\] 
and the theorem follows. 
\end{proof} 

It follows from Theorem \ref{independence} that the Chern character from Theorem \ref{theoremtwo} 
is independent of choice of connection.  We get a corollary:

\begin{corollary} \label{corollaryone} 
 Let $A$ be a smooth  $k$-algebra of finite type where $k$ is a field of characteristic zero.
There exists a ring homomorphism
\[ ch^A: \K_0(A) \rightarrow \H^*_{\DR}(A) .\]
\end{corollary}
\begin{proof} There exists a natural map 
\[ \Omega^p_A \rightarrow
(\Omega^p_A)^{**}=\Hom_A(\wedge^p\Der_k(A),A) \]
hence we get when $\Omega^1_A$ is locally free an isomorphism $i_p:\H^p_{\DR}(A)\cong \H^p(\Der_k(A),A) $. 
Any connection 
\[ \nabla: E\rightarrow E\otimes \Omega^1_A \]
gives rise to a covariant derivation 
\[ \overline{\nabla}:\Der_k(A)\rightarrow \End_k(E) .\] 
One checks that the Chern class defined by
$\overline{\nabla}$ agrees with the one defined by $\nabla$ via $i_p$, and the
claim follows. 
\end{proof} 
The ring homomorphism from Corollary \ref{corollaryone} is the
classical Chern character from Theorem \ref{classicalchern}. 

Note that by functoriality there always exist a diagram
\[ \diagram \K_0(A) \rto \dto^{ch^A} & \K_0(\lg) \dto^{ch^\lg} \\
         \H^*_{\DR}(A) \rto & \H^*(\lg,A) \enddiagram ,\]
but the map $\K_0(A)\rightarrow \K_0(\lg)$ is not surjective in
general: by the example in \cite{maa1}, section 2 the following
holds. Let $k$ be a field of characteristic zero and consider $\O(d)$
on $\P^1_k$. There exist a left $\O_{\P^1}$-linear splitting 
\[ \Pr^1(\O(d))\cong \O(d-1)\oplus \O(d-1) ,\] 
hence the Atiyah-sequence
\[ 0 \rightarrow \Omega^1 \otimes \O(d) \rightarrow \Pr^1(\O(d))
\rightarrow \O(d) \rightarrow 0 \] 
is not left split. It follows that $\O(d)$ does not have a connection. If
we consider the linear Lie-Rinehart algebra $\mathbf{V}_{\O(d)}$ of
$\O(d)$ introduced in section 1 in \cite{maa1}, we see that $\O(d)$ has a $\mathbf{V}_{\O(d)}$-connection.
It follows that the natural map
\[ \K_0(\P^1_k)\rightarrow \K_0( \mathbf{V}_{\O(d)} ) \] 
is not surjective hence $ch^\lg$ is not determined by
$ch^A$ in general. 
Note also that the construction of the Chern-class $ch_n(W,\nabla)$ is valid for any
$S$-algebra A, where $S$ and $A$  are commutative rings. The Chern character
exists  when $S$ is a ring containing the rationals.

\textbf{Acknowledgements}. 
This paper was written in autumn 2003/spring 2004
while the author was a lecturer at the Royal Institute in Stockholm, 
and it is a pleasure to thank Torsten Ekedahl for discussions and remarks on the topics considered.
Thanks also to David Kazhdan and Rolf K{\"a}llstr{\"o}m for valuable comments and inspiring conversations. 
The paper was revised in spring 2005 during a stay at Universite Paris VII
financed by a fellowship of the RTN network HPRN-CT-2002-00287, 
Algebraic K-theory, Linear Algebraic Groups and Related Structures, 
and I would like to thank Max Karoubi for an invitation to Paris.

\end{document}